\documentclass[12pt]{amsart}
\usepackage{amsmath,amssymb,enumerate}
\usepackage{epsf,epsfig,amsfonts,graphicx,color,url}
\usepackage[autostyle]{csquotes}
\usepackage{multirow}

\usepackage{enumitem}
\usepackage[pagewise]{lineno}%\linenumbers
\newcommand{\beq}{\begin{equation}}
\newcommand{\eeq}{\end{equation}}
\usepackage[all]{xy}
\usepackage{filecontents}

\def\R{\mathbb{R}}

\newcommand{\Laa}{\mathfrak{a}}

\newcommand{\Gc}{\mathbb{G}}

\newcommand{\Ag}{\mathbb{A}}
\newcommand{\Ho}{\mathcal{H}}
\newcommand{\Di}{\mathcal{D}}

\numberwithin{equation}{section}

\newtheorem{theorem}{Theorem}

\newtheorem{corollary}{Corollary}[section]

\newtheorem{definition}{Definition}[section]

\renewcommand{\emph}[1]{{\bfseries\itshape{#1}}}
\numberwithin{figure}{section}

\def\sR{subRiemannian } 

\def\ma{metabelian }

\newcommand{\Je}{J^2(\mathbb{R}^2,\mathbb{R})}

\usepackage{hyperref}
\usepackage[%
	capitalize,nameinlink
]{cleveref}
\hypersetup{%
	colorlinks=true,
	linkcolor=blue,
	citecolor=magenta,
	urlcolor=magenta,
}

\begin{document}

\newtheorem*{backgroundtheorem}{Background Theorem}

% the `*' in front gets rid of the numbering;  if I put this above \begin{document}
%formatting gets messed up

\title[Non-integrable subRiemannian geodesic flow on $J^{2}(\R^2,\R)$]{Non-integrable subRiemannian geodesic flow on $J^{2}(\R^2,\R)$}  
\author[A.\ Bravo-Doddoli]{Alejandro\ Bravo-Doddoli} 
\address{Alejandro Bravo: Dept. of Mathematics, UCSC,
1156 High Street, Santa Cruz, CA 95064}
\email{Abravodo@ucsc.edu}
\begin{abstract} 
The space of $2$-jets of a real function of two real variables, denoted by $\Je$, admits the structure of a  metabelian Carnot group, so $\Je$ has a normal abelian sub-group $\Ag$. As any \sR manifold, $\Je$ has an associated Hamiltonian geodesic flow. The Hamiltonian action of $\Ag$ on $T^*\Je$ yields the reduced Hamiltonian $H_{\mu}$ on $T^*\Ho \simeq T^*(\Je/\Ag)$, where $H_{\mu}$ is a two-dimensional Euclidean space. The paper is devoted to proving that reduced Hamiltonian $H_{\mu}$ is non-integrable by meromorphic functions for some values of $\mu$. This result suggests the \sR geodesic flow on $J^{2}(\R^2,\R)$ is not meromorphically integrable.
\end{abstract}

\maketitle

\section{Introduction}

Let $\Je$ be the space of $2$-jets of a real function of two variables. $\Je$ is a Carnot group with step 3 and growth vector $(5,7,8)$. Let $\mathfrak{j}$ be the graded Lie algebra of $\Je$, that is,
$$ \mathfrak{j} = \mathfrak{j}_1 \oplus \mathfrak{j}_2 \oplus \mathfrak{j}_3, \;\; \text{such that} \;\; [\mathfrak{j}_{i},\mathfrak{j}_j] \subseteq \mathfrak{j}_{i+j}\:\;\; \text{and} \;\; \mathfrak{j}_4 = \{ 0 \}. $$
Let $\pi: \Je \to \R^{5}\simeq \mathfrak{j}_1$ be the canonical projection and let $\R^{5}$ be endowed with the Euclidean metric. Consider the \sR metric on  $\Je$ such that $\pi$ is a \sR submersion, see Definition \ref{def:sR-submersion} for the formal definition of a \sR submersion, by construction the \sR structure is left-invariant under the Carnot group multiplication. Like any \sR structure, the cotangent bundle $T^*\Je$ is equipped with a Hamiltonian system whose underlying Hamiltonian $H_{sR}$ is one whose solutions curves are \sR geodesics on $\Je$. This Hamiltonian system is called the \sR geodesic flow on $\Je$.

We say that a group $\Gc$ is \ma if $[\Gc,\Gc]$ is abelian. In \cite{ABD-SY}, we considered the \sR geodesics flow on a general \ma Carnot group $\Gc$. Then, we performed the symplectic reduction of the cotangent bundle $T^*\Gc$ by the Hamiltonian action of the maximal normal abelian sub-group  $\Ag$ containing $[\Gc,\Gc]$, where $\Ag$ acts on $\Gc$ by left multiplication. This action is free and proper, so $\Ho := \Gc/\Ag$ is well defined. Let $\Laa$ be the Lie algebra of $\Ag$ and let $\mu$ be in $\Laa^*$, since $\Ag$ is abelian, the isotropic sub-group of $\Ag_{\mu}:= \{ g \in \Ag: Ad_g^* \mu = \mu\}$ is $\Ag$ and the symplectic reduced space is diffeomorphic to $T^*(\Gc/\Ag) \simeq T^* \Ho$. For more details about the symplectic reduction of the \sR geodesics flow on a \ma Carnot group, see \cite{ABD-SY}, and for the general theory, see \cite{Sym-reduc} or \cite{tudor}.

In the case $\Gc = \Je$, we have $[\mathfrak{j},\mathfrak{j}] = \mathfrak{j}_2\oplus \mathfrak{j}_3$ and the Lie bracket relations in equations \eqref{eq:lie-bra-l2} and \eqref{eq:lie-bra-l3}, see below, show $[\mathfrak{j}_2\oplus \mathfrak{j}_3, \mathfrak{j}_2\oplus \mathfrak{j}_3] = 0$ meaning  $\Je$ is a \ma Carnot group. Following the notation used in \cite{ABD-SY}: $\Ag$ is a 6-dimensional sub-group, whose Lie algebra is framed by $\{E^{\Laa}_{1},E^{\Laa}_{2},E^{\Laa}_{3},E^{\Laa}_{4},E^{\Laa}_{5},E^{\Laa}_6\}$, see equations \ref{eq:lie-bra-l2} and \ref{eq:lie-bra-l3}. Then, $\Ho$ is a 2-dimensional Euclidean space and the reduced Hamiltonian is a two degree of freedom system with polynomial potential, see equation \ref{eq:ham-2-deg}.

The main Theorem of this paper is the following.

\begin{theorem}\label{the:main}
Let $H_{\mu}: T^* \Ho \to \R$ be the reduced Hamiltonian given by the symplectic reduction of \sR geodesic flow on $\Je$ under the action of $\Ag$, where $\mu$ is in $\Laa^*$.  Then, there exists a one parameter family in $\Laa^*$ such that the reduced Hamiltonian $H_{\mu}$ is not meromorphically integrable.
\end{theorem}
Theorem \ref{the:main} suggests the \sR geodesic flow on $J^{2}(\R^2,\R)$ is not meromorphically integrable.

%If $\mu$ in $\Laa^*$ is such that $H_{\mu}$ has an homogeneous polynomial potential of degree $4$.

Examples of Carnot groups with a non-integrable geodesic flow are the following: One is the group of all $4 $ by $4$ lower triangular matrices with 1’s on the diagonal proved by R. Montgomery, M. Shapiro and A. Stolin, see \cite{Shapiro}. Another is the Carnot group with growth vector $(3,6,14)$ showed by I. Bizyaev, A. Borisov, A. Kilin, and I. Mamaev, see \cite{Borisov}. Finally, there is the Carnot group with growth vector $(2,3,5,8)$. Verified by L. V. Lokutsievskiy and Y. L. Sachkov, see \cite{Lokutsievskiy_2018}. 

Kruglikov, B., Vollmer, A. and Lukes-Gerakopoulos, G. made a classification of the integrable geodesic flow on Carnot groups of rank $2$ and low dimension, see \cite{2-rank}.

\section{$\Je$ as a Carnot group}

The  $2$-jet  of a smooth  function $f: \R^2 \to \R$ at a point $(x_0,y_0) \in \R^2$ is its $2$-th order Taylor polynomial at $x_0$.  We will 
encode the  $2$-jet as a   $8$-tuple of real numbers $(j^2 f )|_{(x_0,y_0)}$ as follows: 
\begin{equation*}
(j^2 f )|_{(x_0,y_0)} := \left(x_0,y_0, \frac{\partial^2 f}{\partial x^2},\frac{\partial^2 f}{\partial x \partial y},\frac{\partial^2 f}{\partial y^2} ,\frac{\partial f}{\partial x} ,\frac{\partial f}{\partial  y},f\right)|_{(x_0,y_0)}   \in \R^{8}
\end{equation*}
As $f$ varies over  smooth functions and $(x_0,y_0)$ varies over $\R^2$,  these  $2$-jets
sweep out the   $2$-jet space, denoted   by $\Je$. One can see that $\Je$   is diffeomorphic to $\R^{8}$
and   its   points are coordinatized according to
 $$g = (x,y,u_{2,0},u_{1,1},u_{0,2},u_{1,0},u_{0,1},u) \in \mathbb{R}^{8}.$$
Recall that if $u = f(x,y)$, then $u_{1,0}= \frac{\partial u}{ \partial x}$, $u_{0,1}= \frac{\partial u}{\partial y}$, $u_{2,0}= \frac{ \partial u_{1,0}}{\partial x}$, $u_{1,1} = \frac{\partial u_{1,0}}{\partial y}  = \frac{\partial u_{0,1}}{\partial x}$ and   $u_{0,2} = \frac{\partial u_{0,1}}{\partial y}$. We see that   $\Je$  is endowed with a natural rank 5 distribution $\Di \subset T\Je$
characterized by the following Pfaffian equations
$$ u_{1,0}dx + u_{0,1} dy-du  =  u_{2,0}dx + u_{1,1} dy - du_{1,0} = u_{1,1}dx + u_{0,2} dy - du_{0,1}   =  0. $$

A \sR structure on a manifold consists of a non-integrable distribution $\Di$ together with a smooth inner product $(\cdot , \cdot)_{\Je} $ on $\Di$.  We arrive at the \sR structure 
by observing that $\Di$ is  globally framed by 
 \begin{equation}
 \label{eq:frame}
 \begin{split}
X_1 & = \frac{\partial}{\partial x} + u_{1,0} \frac{\partial}{\partial u} + u_{2,0} \frac{\partial}{\partial u_{1,0}} + u_{1,1} \frac{\partial}{\partial u_{0,1}}, \\
X_2 & = \frac{\partial}{\partial y} + u_{0,1} \frac{\partial}{\partial u} + u_{1,1} \frac{\partial}{\partial u_{1,0}}+ u_{0,2} \frac{\partial}{\partial u_{0,1}}, \\
Y_{1} & = \frac{\partial}{\partial u_{2,0}} , Y_{2}  = \frac{\partial}{\partial u_{1,1}}, Y_{3}  = \frac{\partial}{\partial u_{0,2}}.
\end{split}
\end{equation}  
The Canonical projection $\pi$ is defined by
$$ \pi(g) = g \;\; mod \; [\Je,\Je], $$
and in coordinates is given by $\pi(g) = (x,y,u_{2,0},u_{1,1},u_{0,2})$.  
Now the restrictions of the one-forms  $dx, dy,du_{2,0},du_{1,1},du_{0,2}$ to $\Di$  form a global co-frame for $\Di^*$ which is dual to the frame from equation \eqref{eq:frame}. Let us introduce the formal definition of a \sR submersion.
\begin{definition}\label{def:sR-submersion}
Let $(M,\Di_M,(\cdot , \cdot)_M)$ and $(N,\Di_N,(\cdot , \cdot)_N)$ be two \sR manifolds and let $\phi:M \to N$ a submersion, we consider the case $dim(M) \geq dim(N)$. We say that $\phi$ is a \sR submersion if $(\phi)_* \Di_M = \Di_N$ and $(\phi)^* (\cdot , \cdot)_N = (\cdot , \cdot)_M$.  
\end{definition} 
Therefore \sR  metric on $\Je$ making $\pi$ a \sR submersion is given in coordinates by 
$$(\cdot , \cdot)_{\Je} = (dx^2 +  dy^2 + du_{2,0}^2 +du_{1,1}^2+ du_{0,2}^2)|_{\Di}.$$
An equivalent way to define the \sR metric is to declare the left-invariant vector fields from equation \eqref{eq:frame} orthonormal.
For more details about the jet space as Carnot group, see \cite{CarnotJets}.

Let $\{E_1,E_2,E_1^{\Laa},E_2^{\Laa},E_3^{\Laa}\}$ be the base for first layer $\mathfrak{j}_1$, where $X_i(g) = (L_g)_* E_i$ for $i = 1,2$ and $Y_j(g) = (L_g)_* E_j^{\Laa}$ for $j = 1,2,3$. The frame for $\mathfrak{j}_1$  generates the following Lie algebra:
\begin{equation}\label{eq:lie-bra-l2}
\begin{split}
 E^{\Laa}_{4} &:= [E_1,E^{\Laa}_{1}] = [E_2,E^{\Laa}_{2}]  , \;\; E^{\Laa}_{5} := [E_1,E^{\Laa}_{2}] = [E_2,E^{\Laa}_{3}] ,    \\
\end{split}
\end{equation} 
equations \eqref{eq:lie-bra-l2} define the vector corresponding to the second layer $\mathfrak{j}_2$,
\begin{equation}\label{eq:lie-bra-l3}
\begin{split}
E^{\Laa}_6 &:= [E_1,E^{\Laa}_{4}] = [E_2,E^{\Laa}_{5}] , \\
\end{split}
\end{equation} 
equations \eqref{eq:lie-bra-l3} define the vector corresponding to the third layer $\mathfrak{j}_3$. All the other brackets are zero. The Lie bracket relations in equations \eqref{eq:lie-bra-l2} and \eqref{eq:lie-bra-l3} imply that $\Laa$ is framed by $\{E^{\Laa}_{1},E^{\Laa}_{2},E^{\Laa}_{3}$, $E^{\Laa}_{4},E^{\Laa}_{5},E^{\Laa}_6\}$.  Let $\Ho$ be the 2-dimensional Euclidean space defined by quotient $\Je/\Ag$. Since $[E_1,E_2] = 0$, we can think $\Ho$ as a sub-group  of $\Je$ such that $\Je \simeq  \Ag \rtimes \Ho$.

%%
%%In \cite{Sym-reduc}, we introduced the following definition,  we say that $\Gc$ \ma group is a $n$-abelian extension if $\Ho := \Gc /\Ag$ is a $n$-dimensional Euclidean space. Then, by definition, $\Je$ is a $2$-abelian extension, and Theorem \ref{the:main} contributes to the classification of $2$-abelian extension Carnot Groups with the non-integrable geodesic flow.

\subsection{The exponential coordinates of the second kind}

The jet space $\Je$ has a natural definition using the coordinates $x$, $y$, and $u$'s; however, these coordinates do not easily show the symmetries of the system. The exponential coordinates of the second kind exhibit the symmetries:

We recall that the exponential map $\exp:\mathfrak{j}\to \Je$ is a global diffeomorphism, this allow us to endow $\Je$ with coordinates $(x,y,\theta_{1},\theta_{2},\theta_{3},\theta_{4},\theta_{5},\theta_6)$ in the following way: a point $g$ in $\Je$ is given by
\begin{equation*}
g := \prod_{i=1}^{6} \exp(\theta_i E^{\Laa}_{i})*\exp(y E_{2})*\exp(x E_{1}).
\end{equation*}

Then the horizontal left-invariant vector fields are given by  
\begin{equation}\label{eq:vec-firts-layer-h}
\begin{split}
X_1 & := \frac{\partial}{\partial x} , \qquad X_2 := \frac{\partial}{\partial y},\\
\end{split}
\end{equation} 
and
\begin{equation}\label{eq-lef-in-can}
\begin{split}
 Y_{1} & :=  \frac{\partial}{\partial \theta_{1}} + x \frac{\partial}{\partial \theta_{4}}  + \frac{x^2}{2!}\frac{\partial}{\partial \theta_6},  \\
Y_{2} & :=  \frac{\partial}{\partial \theta_{2}}  + y \frac{\partial}{\partial \theta_4} + x \frac{\partial}{\partial \theta_{5}} + xy \frac{\partial}{\partial \theta_6}, \\
Y_{3} & :=  \frac{\partial}{\partial \theta_{3}} + y \frac{\partial}{\partial \theta_{5}} + \frac{y^2}{2!}  \frac{\partial}{\partial \theta_6}. \\
\end{split}
\end{equation} 
The left-invariant vector fields from equation \eqref{eq:vec-firts-layer-h} and \eqref{eq-lef-in-can} just depend on the independent variables $x$ and $y$. All the \ma Carnot groups have this property, which is the heart of the symplectic reduction. For more details, see \cite{ABD-SY}. 

The change from the coordinates $(x,y,u_{2,0},u_{1,1},u_{0,2},u_{1,0},u_{0,1},u)$ to the exponential coordinates of the second kind  $(x,y,\theta_{1},\theta_2,\theta_3,\theta_4,\theta_5,\theta_6)$ is the following
\begin{equation*}
\begin{pmatrix}
x \\
y \\
\theta_{1} \\
\theta_2 \\ 
\theta_3 \\
\theta_4 \\
\theta_5 \\
\theta_6 \\
\end{pmatrix}
=
\begin{pmatrix}
x \\
y \\
u_{2,0} \\ 
u_{1,1} \\
u_{0,2} \\
x u_{2,0} + yu_{1,1} - u \\
x u_{1,1} + yu_{0,2} - u \\
\frac{x^2}{2} + xy u_{1,1} + \frac{y^2}{2} u_{0,2} - x u_{1,0} - yu_{0,1} + u \\
\end{pmatrix}
\end{equation*}

\section{Geodesic flow on $\Je$}

Let us consider the traditional coordinates on $T^*\Je$, that is, $(p,g)$ where $p:= (p_x,p_y,p_1,p_2,p_3,p_4,p_5,p_6)$ are the momentum associated with exponential coordinates of the second type, see \cite{Arnold} and \cite{Landau} for more details about the traditional coordinates. Let $ P_{X_1},$ $P_{X_2}$, $P_{Y_{1}}$, $P_{Y_2}$ and $P_{Y_{3}}$ be the momentum functions associated with the left-invariant vector fields on the first layer $\mathfrak{j}_1$ are given by
\begin{equation}
\begin{split}
P_{X_1} & = p_x, \qquad P_{X_2} = p_y\qquad
 Y_{1} =  p_{1} + x p_{4}  + \frac{x^2}{2!}p_6,  \\
Y_{2} & = p_{2}  + y p_4 + x p_{5} + xy p_6, \qquad Y_{3}  =  p_{3} + y p_{5} + \frac{y^2}{2!}  p_6, \\
\end{split}
\end{equation}
see \cite{tour}, or \cite{agrachev} for more details about the momentum functions. Then, the Hamiltonian governing the \sR geodesic flow on $\Je$ is
\begin{equation}
H_{sR} := \frac{1}{2}(P_{X_1}^2+P_{X_2}^2+P_{Y_1}^2+P_{Y_2}^2+P_{Y_3}^2).
\end{equation}
See \cite{tour}, or \cite{agrachev} for more details about the definition of $H_{sR}$.

The Hamiltonian function $H_{sR}$ does not depend on the coordinates $\theta_{1}$, $\theta_{2}$, $\theta_{3}$, $\theta_{4}$, $\theta_{5}$ and $\theta_6$, so they are cyclic coordinates, in other words, $p_1$, $p_2$, $p_3$, $p_4$, $p_5$ and $p_6$ are constants of motion, see \cite{Landau} or \cite{Arnold} for more details about the cyclic coordinates. Moreover, since $H_{sR}$ is invariant under the action of $\Ag$, these constants of motion correspond to the momentum map $J:T^*\Je \to \Laa^*$ given by
$$ J(p,g) = \mu := (a_1,a_2,a_3,a_4,a_5,a_6) \;\; \text{where}\;\;  p_i = a_i, \;\; 1 \leq i \leq 6.$$ 
See \cite{tour} or \cite{tudor} for the formal definition of the momentum map. See \cite{ABD-SY} for the construction of the momentum map in the context of \ma Carnot group.

\subsection{The reduced Hamiltonian }\label{sub-sect:reduce}
By the symplectic theory, the reduced space is diffeomorphic to $T^*(G/A) \simeq T^* \Ho$, and the reduced Hamiltonian is a two-degree-of-freedom system with a polynomial potential of degree four in the variables $x$ and $y$, and depending on the parameters $\mu  :=  (a_{1},a_{2},a_{3},a_{4},a_{5},a_6)$ in $\Laa^*$, given by
\begin{equation}\label{eq:ham-2-deg}
H_{\mu}(p_x,p_y,x,y) := \frac{1}{2} \left( p_x^2 + p_2^2 + \phi_{\mu}(x,y) \right),
\end{equation}
where $\phi_{\mu}(x,y)$ is the following potential
\begin{equation}\label{eq:pot}
 (a_{1} + a_{4} x + \frac{x^2}{2!} a_6)^2 + ( a_{2} + a_{5} x + a_{4} y + a_6 x y )^2 + (a_{3} + a_{5} y + a_6 \frac{y^2}{2!} )^2.
\end{equation}
Let $\pi_{\Ag}:\Je \to \Ho$ be the canonical projection given by 
$$\pi_{\Ag}(g)= (x,y).$$
Let $\Pi_{\Ag}: T^*G \to T^*\Ho$ be co-lift projection associated to $\pi_{\Ag}$, that is,
$$\Pi_{\Ag}(p,g) = (p_x,p_y,x,y) .$$ 
Then symplectic reduction implies   
$$H_{sR}|_{J^{-1}(\mu)} = H_{\mu} \circ \Pi_{\Ag}.$$

\subsection{Background Theorem}

Here we introduce the \textbf{ Background Theorem}, which provides a complete classification of the Yang-Mills Hamiltonian system by S. Shi and W. Li, in \cite{Yang-Mills}.
\begin{backgroundtheorem}
\label{thm:back} 
Let $H$ be the Hamiltonian system given by
\begin{equation}
H = \frac{1}{2c}(p_x^2+p_y^2) + \frac{1}{2c}(ax^2+by^2) + \frac{1}{4c^2}(cx^2+dy^2+2ex^2y^2),
\end{equation}
where $a$, $b$, $c \neq 0$, $d$ and $e$ are in $\R$. Then $H$ is meromorphically integrable in the Liouvillian sense (i.e., the existence of an additional meromorphic integral) if and only if one of the following conditions hold:
\begin{enumerate}[label=(\Alph*)]
\item $e=0$,
\item $c = d = e$,
\item $a = b$ and $e = 3c = ed$,
\item $b = 4a$, $e = 3c$ and $d = 8c$,
\item $b = 4a$, $e = 6c$ and $d =16c$,
\item $b = 4a$, $e = 3d$ and $c = 8d$,
\item $b = 4a$, $e = 6d$ and $c = 16d$.
\end{enumerate}
\end{backgroundtheorem}

A consequence of the \textbf{ Background Theorem} is the following.
\begin{corollary}\label{cor}
Let $H$ be the Hamiltonian given by
\begin{equation}\label{eq:cor}
H = \frac{1}{2}(p_x^2+p_y^2) + \frac{1}{4c^2} ( c x^4 + 4 c x^2y^2 + c y^4).
\end{equation} 
Then $H$ is not meromorphically integrable in the Liouvillian sense.
\end{corollary}
\begin{proof}
Following the notation from the \textbf{Background Theorem}, we have that $c = d$, $e = 2c$, and $a = b = 0$, so $H$ is not meromorphically integrable. 
\end{proof}
An alternative proof of corollary \ref{cor} is given by
Maciejewski, A. J.  and Przybylska, M., in \cite{Darboux-points}.

\subsection{Proof of Theorem \ref{the:main}}

Now we are ready to prove Theorem \ref{the:main}.

\begin{proof}
If $\mu  =   \left(0,0,0,0,0,0,0,\sqrt{\frac{2}{c}}\right)$ with $c$ in $(0,\infty)$, then equation \eqref{eq:pot} implies the potential $\phi_{\mu}(x,y)$ is $ \frac{1}{4c^2} (c x^4 + 4c x^2 y^2 +c y^4)$. Let $H_{\mu}$ be given by equation \eqref{eq:ham-2-deg}, then  $H_{\mu}$ is equal to the Hamiltonian given by equation \eqref{eq:cor}, so by Corollary \ref{cor} $H_{\mu}$ is not integrable by meromorphic functions.  
\end{proof}

\begin{appendix}

%\section{\sR submersion}

%
%\section{The $\Laa^*$ value one-form $\alpha_{\Je}$ }
%
%In \cite{ABD-SY}, we showed that the mathematical object relating the \sR geodesic flow on $\Je$ and the reduced Hamiltonian on $T^*\Ho$ is $\Laa^*$ value one-form $\alpha_{\Je}$ on  $\mathfrak{j}_1 \simeq \R^5$ given by
%\begin{equation}\label{eq:alp-Je}
%\begin{split}
% \alpha_{\Je} & = d\theta_{1} \otimes (e_{1} + x e_{4}  + \frac{x^2}{2!} e_6) \\
%  & + d\theta_{2} \otimes ( e_{2} + x e_{5}  + y e_{4}  +  x y e_6) \\
%  & + d\theta_{3} \otimes (e_{3} + y e_{5}  + \frac{y^2}{2!} e_6 ).  \\
%\end{split} 
%\end{equation}

\end{appendix}

\nocite{*} % to test all bib entrys
\bibliographystyle{plain}
\bibliography{bibli} %

\begin{thebibliography}{10}

\bibitem{agrachev}
Barilari~D. Agrachev, A. and U.~Boscain.
\newblock {\em A Comprehensive Introduction to Sub-Riemannian Geometry}.
\newblock Cambridge University Press, 2019.

\bibitem{Arnold}
V.~I. Arnol'd.
\newblock {\em Mathematical methods of classical mechanics, 2nd ed., New
  York:}.
\newblock Springer, 1989.

\bibitem{Borisov}
I.~Bizyaev, A.~Borisov, A.~Kilin, and I.~Mamaev.
\newblock Integrability and nonintegrability of sub-riemannian geodesic flows
  on carnot groups.
\newblock {\em Regular and Chaotic Dynamics}, 21:759--774, 11 2016.

\bibitem{2-rank}
Boris Kruglikov, Andreas Vollmer, and Georgios Lukes-Gerakopoulos.
\newblock On integrability of certain rank 2 sub-riemannian structures.
\newblock {\em Regular and Chaotic Dynamics}, 2015.

\bibitem{Landau}
L.D. Landau and E.M. Lifshitz.
\newblock {\em Mechanics third edition: Volume 1 of course of theoretical
  physics}.
\newblock 1976.

\bibitem{Lokutsievskiy_2018}
L.~V. Lokutsievskiy and Yu.~L. Sachkov.
\newblock Liouville integrability of sub-riemannian problems on carnot groups
  of step 4 or greater.
\newblock {\em Sbornik: Mathematics}, 209(5):672, 2018.

\bibitem{Darboux-points}
A.~J. Maciejewski and M.~Przybylska.
\newblock Darboux points and integrability of hamiltonian systems with
  homogeneous polynomial potential.
\newblock {\em Journal of Mathematical Physics}, 46(6):062901, 2005.

\bibitem{Sym-reduc}
Jerrold Marsden and Alan Weinstein.
\newblock Reduction of symplectic manifolds with symmetry.
\newblock {\em Reports on Mathematical Physics}, 5, 1974.

\bibitem{tour}
R.~Montgomery.
\newblock {\em A Tour of Subriemannian Geometries, Their Geodesics and
  Applications}.
\newblock Number~91. American Mathematical Soc.

\bibitem{Shapiro}
R.~Montgomery, M.~Shapiro, and A.~Stolin.
\newblock Chaotic geodesics in carnot groups.
\newblock preprint on webpage at \url{https://arxiv.org/abs/dg-ga/9704013},
  1997.

\bibitem{tudor}
Juan-Pablo Ortega and Tudor Ratiu.
\newblock {\em Momentum Maps and Hamiltonian Reduction}.
\newblock Progress in Mathematics, 2004.

\bibitem{ABD-SY}
N.~Paddeu, A.~Bravo-Doddoli, and E.~Le Donne.
\newblock Sympletic reduction of the sub-riemannian geodesic flow on
  meta-abelian carnot groups.
\newblock preprint on webpage at \url{https://arxiv.org/abs/2211.05846}, 2022.

\bibitem{Yang-Mills}
S.~Shi and W.~Li.
\newblock Non-integrability of generalized yang-mills hamiltonian system.
\newblock {\em Discrete and Continuous Dynamical Systems}, 33(4):1645--1655,
  2013.

\bibitem{CarnotJets}
B.~Warhurst.
\newblock Jet spaces as nonrigid carnot groups.
\newblock {\em Journal of Lie Theory}, 15(1):341--356, 2005.

\end{thebibliography}

\end{document}